\documentclass[12pt,dvips]{article}
\usepackage{amsmath}
\usepackage{amssymb}
\usepackage{amsthm}
\usepackage{latexsym}
\usepackage[all,line,rotate]{xy}
\usepackage{tabls}

\DeclareMathOperator{\interior}{int}

\DeclareMathOperator{\cone}{cone}

\def\vec#1{\mathchoice{\mbox{\boldmath$\displaystyle\bf#1$}}
{\mbox{\boldmath$\textstyle\bf#1$}}
{\mbox{\boldmath$\scriptstyle\bf#1$}}
{\mbox{\boldmath$\scriptscriptstyle\bf#1$}}}
\newdimen\intercol

\begin{document}
\theoremstyle{plain}
\newtheorem{theorem}{Theorem}
\newtheorem{prop}[theorem]{Proposition}
\newtheorem{cor}[theorem]{Corollary}
\newtheorem{lemma}[theorem]{Lemma}
\newtheorem{question}[theorem]{Question}
\newtheorem{conj}[theorem]{Conjecture}
\newtheorem{assumption}[theorem]{Assumption}

\theoremstyle{definition}
\newtheorem{definition}[theorem]{Definition}
\newtheorem{notation}[theorem]{Notation}
\newtheorem{condition}[theorem]{Condition}
\newtheorem{example}[theorem]{Example}
\newtheorem{introduction}[theorem]{Introduction}

\numberwithin{theorem}{section}

\makeatletter                          % This sequence of commands will
\let\c@equation\c@theorem              % incorporate equation numbering
\makeatother                           % into the theorem numbering scheme
\renewcommand{\theequation}{\arabic{section}.\arabic{equation}}

\makeatletter                          % This sequence of commands will
\let\c@figure\c@theorem              % incorporate figure numbering
\let\c@table\c@theorem              % incorporate table numbering
\makeatother                           % into the theorem numbering scheme
\renewcommand{\thefigure}{\arabic{section}.\arabic{figure}}
\renewcommand{\thetable}{\arabic{section}.\arabic{table}}

\newcommand{\todo}[1]{\vspace{5 mm}\par \noindent
   \marginpar{\textsc{ToDo}}\framebox{\begin{minipage}[c]{0.95 \textwidth}
   \tt #1 \end{minipage}}\vspace{5 mm}\par}

\providecommand{\abs}[1]{\left\lvert#1\right\rvert}
\providecommand{\norm}[1]{\left\lVert#1\right\rVert}
\providecommand{\floor}[1]{\left\lfloor#1\right\rfloor}
\providecommand{\Z}{\mathbb{Z}} \providecommand{\R}{\mathbb{R}}
\providecommand{\N}{\mathbb{N}} \providecommand{\C}{\mathbb{C}}
\providecommand{\Q}{{\mathbb{Q}}} \providecommand{\x}{\mathbf{x}}
\providecommand{\y}{\mathbf{y}} \providecommand{\z}{\mathbf{z}}
\providecommand{\boxend}{\hspace{\stretch{1}}$\Box$\\ \ \\}
\long\def\symbolfootnote[#1]#2{\begingroup%
\def\thefootnote{\fnsymbol{footnote}}\footnote[#1]{#2}\endgroup}
\renewcommand{\mod}[1]{\,(\text{mod }#1)}

\providecommand{\D}{\mathcal{D}}

\title{Counting with rational generating functions}
\author{Sven Verdoolaege and Kevin Woods}
\maketitle

\begin{abstract}
We examine two different ways of encoding a counting function, as a
rational generating function and explicitly as a function (defined
piecewise using the greatest integer function).  We prove that, if the degree and number of input variables of the (quasi-polynomial) function are fixed,
there is a polynomial time algorithm which converts
between the two representations.  Examples of such counting
functions include Ehrhart quasi-polynomials, vector partition
functions, integer points in parametric polytopes, and projections
of the integer points in parametric polytopes.  For this last
example, this algorithm provides the first known way to
compute the explicit function in polynomial time.
We rely heavily on results of
Barvinok, and also of Verdoolaege, Seghir, Beyls, et al.
\end{abstract}

\section{Introduction}

We are interested in a wide variety of functions of the form
$$
c:\Z^{n}\rightarrow\Q
.
$$
Most examples, including
Ehrhart quasi-polynomials and vector partition functions, will count some combinatorial object.
The function $c(\vec s)$
can be encoded in at least two different ways: either as
an explicit function or as a generating function
\[f(\x) = \sum_{\vec s=(s_1,\ldots,s_n)\in\Z^n}c(\vec s)x_1^{s_1}\cdots x_n^{s_n}=\sum_{\vec s \in
\Z^{n}} c(\vec s) \x^{\vec s}.\]

\begin{example}
\label{BabyExample}
Consider the generating function
\[
f(x)=\frac{1}{1-x^2}=1+x^2+x^4+\cdots = \sum_{s\in\Z}c(s) x^{s} .\]
The corresponding function can be represented explicitly as
\[c(s)=\left\{%
\begin{array}{ll}
    0, & \hbox{if $s<0$} \\
    0, & \hbox{if $s\ge 0$ and $s$ odd} \\
    1, & \hbox{if $s\ge 0$ and $s$ even}.
\end{array}%
\right.
\]
\end{example}\boxend

%Furthermore,
%building on earlier results from \cite{BW03} and using this
%polynomial conversion, we will show that we can compute an explicit
%function representation for $c(\vec s)$ in polynomial time (for
%fixed dimensions) provided $P$ is a polytope.

Mathematicians often encode a function as a rational
generating function, such as $f(x)=\frac{1}{1-x^2}$ in Example
\ref{BabyExample}, which is a compact representation of a (possibly
infinite) Laurent power series $\sum_{\vec s\in\Z^n}c(\vec
s)\x^{\vec s},$ where $c(\vec s)\in\Z^n$.  This has the advantage
that we may apply many computational tools to manipulate our
rational generating function and obtain information from it (see,
for example, \cite{BW03}). An explicit function representation for
$c(\vec s)$, on the other hand, has the advantage of being easily
evaluated for a particular value of $\vec s$. Such a representation
is therefore preferred in the compiler community (see, for example,
\cite{Verdoolaege2007parametric}).

We will show that these ways of representing a function are
``the same,'' in the sense that one can convert between the rational
function and explicit function representations in polynomial time (if the degree and number of variables of the function is fixed). Let us be more precise about the specific representations we will
use for generating functions and explicit functions.
\begin{definition}
\label{def:generating}
By a \emph{rational
generating function} $f(\x)$, we will mean a function given to us in
the form
\begin{equation} \label{GFForm}
f(\x)=\sum_{i\in I}\alpha_i \frac{\x^{\vec p_i}}
     {(1-\x^{\vec b_{i1}})(1-\x^{\vec b_{i2}})\cdots(1-\x^{\vec b_{ik_i}})},
\end{equation}
where $\x\in \C^n$, $I$ is a finite set, $\alpha_i\in\Q$,
$\vec p_i\in\Z^n$, and $\vec b_{ij}\in\Z^n\setminus \{\vec 0\}$.
\end{definition}
\begin{definition}
A {\em step-polynomial\/} $g : \Z^n \to \Q$ is a function written in
the form
\[
g(\vec s) = \sum_{j=1}^{m}\alpha_j\prod_{k=1}^{d_j} \floor{\langle
\vec a_{jk},\vec s\rangle+b_{jk}},
\]
where $\alpha_j\in\Q$, $\vec a_{jk}\in\Q^n$, $b_{jk}\in\Q$,
$\langle \cdot,\cdot\rangle$ is the standard inner product, and
$\floor{\cdot}$ is the greatest integer function. We say that the
\emph{degree} of $g(\vec s)$ is $\max_{j}\{\,d_{j}\,\}$.

A {\em piecewise step-polynomial\/} $c:\Z^n\rightarrow\Q$ is a
collection of polyhedra $Q_i$ (which may not all be full dimensional) together with
corresponding functions $g_i:Q_i \cap \Z^n \rightarrow \Q$ such
that\begin{enumerate}
    \item the $\interior(Q_i)$ partition $\Q^n$ (where $\interior(Q)$ is the relative interior of $Q$ in the affine space it lies in)
    \item $c(\vec s)=g_i(\vec s)$, for $\vec s\in \interior(Q_i) \cap \Z^n$, and
    \item each $g_i$ is a step-polynomial.
\end{enumerate}
We say that the degree of $c(\vec s)$ is $\max_i \deg g_i$.  Working with the relative interiors of the polyhedra allows us not to worry about the value of the function at the intersection of two polyhedra.
\end{definition}
For example, the explicit function $c(s)$ in Example~\ref{BabyExample}
can be written as the piecewise step-polynomial
$$
c(s) =
\begin{cases}
1 + \left\lfloor \frac s 2 \right\rfloor - \left\lfloor \frac {s+1}
2 \right\rfloor, & \hbox{if $s > 0$}\\
1,& \hbox{if s=0}\\
0, & \hbox{if $s<0$}
.
\end{cases}
$$

We must be careful when speaking of a correspondence between a
rational generating function and a piecewise step-polynomial,
because a generating function may have different Laurent power
series expansions which converge on different regions of $\C^n$. For
example, if $f(x)=\frac{1}{1-x}$ then
\[1+x+x^2+x^3+\cdots \text{ and } -x^{-1}-x^{-2}-x^{-3}-\cdots\]
are Laurent power series expansions convergent for $\norm{x}<1$ and
$\norm{x}>1$, respectively.

We state the main theorem, and then provide several examples of
rational generating functions and piecewise step-polynomials.

\begin{theorem}
\label{ThmEquivalence} Fix $n$ and $k$. There is a polynomial time
algorithm which, given a rational generating function $f(\x)$ in the
form (\ref{GFForm}) with $n$ variables and each $k_i\le k$ and given
$\vec l\in\Z^n$ such that $\langle \vec l, \vec b_{ij}\rangle \ne
0$ for all $i$ and $j$, computes the piecewise step-polynomial
$c:\Z^n\rightarrow\Q$ with degree at most $k$ such that
\[f(\x)=\sum_{\vec s\in\Z^n}c(\vec s)\x^{\vec s}\] is the Laurent power series
expansion of $f(\x)$ convergent on a neighborhood of
$\mathbf{e}^{\vec l}=(e^{l_1},e^{l_2},\ldots,e^{l_n})$,
with $e$ the base of the natural logarithmic function.

Conversely, there is a polynomial time algorithm which, given a
piecewise step-polynomial $c:\Z^n\rightarrow\Q$ of degree at most
$k$ such that $f(\x)=\sum_{\vec s\in\Z^n}c(\vec s)\x^{\vec s}$
converges on some nonempty open subset of $\C^n$, computes the
rational generating function $f(\x)$ in the form (\ref{GFForm}) with
$k_i\le k$.
\end{theorem}

\ \\

The proof of the first half of this theorem
 will use several ideas from \cite{BP99}. Section 3 will be devoted to
the proof of the theorem, after we lay the groundwork in Section 2.
Note that applying the theorem twice (in one direction and then the other) will in general
not result in the exact same representation of the rational generating function
or piecewise step-polynomial.
We are unaware of any canonical form for either rational generating functions
or piecewise step-polynomials that can be computed in polynomial time.

As there may be many functions with the same generating function
representation (convergent on different neighborhoods), we need
to find an appropriate $\vec l$ value when we want to convert
a given rational generating function to an explicit representation.
If we know that the function $c(\vec s)$ is only nonzero for $\vec s$ in
some polyhedron $Q$ such that $Q$ does not contain any straight
lines, then we may take any $\vec l$ such that $\langle \vec
l,\vec b_{ij}\rangle\ne 0$ for all $i,j$ and such that
\[Q\cap\{\,\vec x\in\Q^n \mid \langle \vec l,\vec x \rangle \ge
0\,\}\] is bounded. Such an $\vec l$ will give us the desired Laurent
power series expansion $\sum_{\vec s}c(\vec s)\,\x^{\vec s}$. In
Example \ref{BabyExample}, we could take $l=-1$.

\begin{example}
Let $P\subset\Q^d$ be a rational polytope, and let
\[c_P(s)=\#(sP\cap\Z^d),\] where $sP$ is $P$ dilated by a factor of
$s$.
\end{example} Then Ehrhart proved \cite{Ehrhart62} that $c_P(s)$
is a quasi-polynomial, that is, there is a $\D\in\Z_+$ and
polynomial functions $g_0(s),g_1(s),\ldots,g_{\D-1}(s)$ such that
\[c_P (s) = g_j(s)\text{ for }s \equiv j \mod{\D}.\]
The generating function $\sum_{s=0}^{\infty} c_P(s)x^s$ can be
computed in polynomial time, and this has been implemented in LattE
(see \cite{DHT03}).  Computing some explicit function representation
of $c_P(s)$ in worst-case exponential time
has been implemented in \cite{CL98} and computing
$c_P(s)$ as a piecewise step-polynomial in polynomial time
had been implemented in \cite{Verdoolaege2007parametric}.

\begin{example}
In particular, let $P\subset \Q^2$ be
$\left[0,\frac{1}{2}\right]\times \left[0,\frac{1}{2}\right]$.
\end{example}
Then
\[c_P(s)=\floor{\frac{1}{2}s+1}^2, \text{ for $s\ge 0$},\]
and we have that
\[\sum_{s=0}^{\infty}
c_P(s)x^s=\frac{2}{(1-x)(1-x^2)^2}-\frac{1}{(1-x)(1-x^2)},\]
which can be verified by hand.
\boxend

\begin{example}
\label{ex:VPFDef} Given
$\vec a_1,\vec a_2,\ldots,\vec a_d\in\N^n$, let
$c:\Z^n\rightarrow \Z$ be the \emph{vector partition function},
defined by
\[c(\vec s)=\#\left\{\,
\vec\lambda=(\lambda_1,\lambda_2,\ldots,\lambda_d)\in\N^d \mid
\vec
s=\lambda_1\vec a_1+\lambda_2\vec a_2+\cdots+\lambda_d\vec a_d\,\right\},\]
i.e., the number of ways an integer vector $\vec s$ can be written
as a nonnegative combination of the $\vec a_i$.
\end{example}
Then the generating function representation of $c(\vec s)$ is very
simple:
\[f(\x)=\frac{1}{(1-\x^{\vec a_1})(1-\x^{\vec a_2})\cdots(1-\x^{\vec a_d})}.\]
The piecewise step-polynomial representation of $c(\vec s)$ can also
be computed in polynomial time (see Corollary \ref{CorVPF} or
\cite{Verdoolaege2007parametric}).
        Beck~\cite{Beck2004fractions} describes a general technique
        for computing vector partition functions, based on partial fraction expansions of $f(\x)$.
        He does not provide a complexity analysis, but
        standard techniques for computing partial
        fractions~\cite{Henrici1974}
        are exponential, even for fixed dimensions.

\begin{example}
In particular, consider the number of ways to partition an integer $s$ into 2's and 5's, i.e.,
$a_1=2$ and $a_2=5$.
Then the generating function representation is
\[f(x)=\frac{1}{(1-x^2)(1-x^5)},\]
and
\[c(s)=\left\{%
\begin{array}{ll}
    0, & \hbox{if $s<0$} \\
    \floor{\frac{1}{2}s+1}+\floor{-\frac{2}{5}s}, & \hbox{if $s\ge 0$},\\
\end{array}%
\right.
\]
which, again, can be verified by hand.
\boxend
\end{example}

Both Ehrhart quasi-polynomials and vector partition functions are
special cases of counting integer points in \emph{parametric
polytopes}. In general, we let $P\subset\Q^{n}\times\Q^{d}$ be a
rational polyhedron such that, for all $\vec s\in\Q^{n}$, the set
$P_{\vec s} = \{\vec t\in\Q^{d} \mid  (\vec s, \vec t)\in P\}$
is bounded, and we define the function $c:\Z^{n}\rightarrow \Z$ by
\begin{equation}
\label{eq:parametric} c(\vec s) = \# (P_{\vec s} \cap \Z^{d}) =
\#\left\{\vec t\in\Z^{d} \mid  (\vec s,\vec t)\in P\right\}.
\end{equation}
We call $P$ a parametric polytope, because, if $P=\big\{(\vec s,\vec
t)\in\Q^{n}\times\Q^{d} \mid A\vec s+B\vec t\le \vec c\big\}$ for
some matrices $A\in\Z^{m\times n}$, $B\in\Z^{m\times d}$ and vector
$\vec c\in\Z^m$, then
\[P_{\vec s}=\left\{\vec t\in\Q^{d} \mid B\vec t\le \vec c-A\vec s\right\},\]
so as $\vec s$ varies, the polytope $P_{\vec s}$ varies by changing
the right hand sides of its defining inequalities.

Both a piecewise step-polynomial representation for $c(\vec s)$ and
its generating function, $\sum_{\vec s}c(\vec s)\x^{\vec s}$, can be
computed in polynomial time, as the following two propositions
state.
\begin{prop}[\cite{Verdoolaege2007parametric}]
\label{p:explicit} Fix $n$ and $d$.  There is a polynomial time
algorithm which, given a parametric polytope $P\subset
\Q^{n}\times\Q^{d}$, computes the piecewise step-polynomial
\[c(\vec s)=\# (P_{\vec s} \cap \Z^{d})\]
with degree at most $d$.
\end{prop}

\ \\

\begin{prop}
\label{p:generating} Fix $n$ and $d$.  There is a polynomial time
algorithm which, given a parametric polytope $P\subset
\Q^{n}\times\Q^{d}$ such that
\[f(\x)=\sum_{\vec s\in\Z^{n}} c(\vec s)\x^{\vec s}\]
converges on some nonempty open subset of $\C^{n}$, computes $f(\x)$
as a rational generating function of the form (\ref{GFForm}) with
the $k_i$ at most $d$.
\end{prop}

\ \\

In Section~\ref{s:explicit} we recall the key ideas of the
proof of Proposition~\ref{p:explicit} from \cite{Verdoolaege2007parametric},
drawing heavily from the ideas in \cite{BP99}.
Proposition~\ref{p:generating} is an immediate consequence of
\cite[Theorem 4.4]{BP99} and the monomial substitution from
\cite{BW03} and will be proved in Section~\ref{s:equivalence}.
%A general overview of these two propositions is shown in
%Figure~\ref{f:overview}.

%\begin{figure}
%\begin{description}
%\item[Input:] A parametric polytope
%$P_{\vec s} = \{\vec t\in\Q^{d} \mid  (\vec s, \vec t)\in P\}$
%with $P\subset\Q^{n}\times\Q^{d}$ a rational polyhedron
%\item[Output:] A piecewise step-polynomial
%$c(\vec s) = \# (P_{\vec s} \cap \Z^{d}$ or
%a rational generating function
%$f(\x)=\sum_{\vec s\in\Z^{n}} c(\vec s)\x^{\vec s}$
%\end{description}
%\begin{itemize}
%\item Computation of $c(\vec s)$ (Proposition~\ref{p:explicit},
%    page~\pageref{p:explicit:proof})
%\begin{itemize}
%\item Compute chamber decomposition (Proposition~\ref{p:chambers},
%    page~\pageref{p:chambers})
%\item For each chamber $C_i$ in decomposition
%\begin{itemize}
%\item Compute $f(P_{\vec s}\cap\Z^d;\x)$ for $\vec s\in C_i$
%(Proposition~\ref{p:core}, page~\pageref{p:core})
%\item Compute $c_i(\vec s):=f(P_{\vec s}\cap\Z^d;\vec 1)$
%(Lemma~\ref{l:specialization}, page~\pageref{l:specialization})
%\end{itemize}
%\end{itemize}
%\item Computation of $f(\x)$ (Proposition~\ref{p:generating},
%    page~\pageref{p:generating:proof})
%\begin{itemize}
%\item Compute $f(P\cap\Z^{n+d};\x,\y)$
%(Proposition~\ref{p:core}, page~\pageref{p:core})
%\item Compute $f(\x) = f\big(P\cap\Z^{n+d};\x, \vec 1\big)$
%(Lemma~\ref{l:specialization}, page~\pageref{l:specialization})
%\end{itemize}
%\end{itemize}
%\caption{General overview of Propositions~\ref{p:explicit}
%and~\ref{p:generating}.}
%\label{f:overview}
%\end{figure}

We may also look at \emph{projections} of the integer points in a
parametric polytope.  Let $P\subset\Q^{n}\times\Q^{d}\times\Q^{m}$
be a rational polytope, and define the function $c:\Z^{n}\rightarrow
\Z$ by
\[c(\vec s)=\#\left\{\vec t\in\Z^{d} \mid \exists \vec u\in\Z^{m}:\
(\vec s,\vec t,\vec u)\in P\right\}.\] If $P_{\vec s}=\left\{(\vec
t,\vec u)\in\Q^d\times\Q^m \mid (\vec s,\vec t,\vec u)\in P\right\}$
and the projection $\pi:\Q^d\times\Q^m\rightarrow\Q^d$ is defined by
$\pi(\vec t, \vec u)=\vec t$, then \[c(\vec s)=\#\big(\pi(P_{\vec
s}\cap\Z^{d+m})\big).\] It follows from \cite{BW03} that the
generating function, $\sum_{\vec s}c(\vec s)\x^{\vec s}$, can be
computed in polynomial time (for fixed $n$, $d$, and $m$).
Therefore, we have as a corollary to Theorem \ref{ThmEquivalence}
that the piecewise step-polynomial can be computed in polynomial
time.

\begin{cor}
\label{CorParametricProjection} Let $n$, $d$, and $m$ be fixed.
There is a polynomial time algorithm
which, given a polytope $P\subset\Q^{n}\times\Q^{d}\times\Q^{m}$,
computes the piecewise step-polynomial
\[c(\vec s)=\#\left\{\vec t\in\Z^{d} \mid \exists \vec u\in\Z^{m}:\
(\vec s,\vec t,\vec u)\in P\right\}\] with degree at most $n+d+m$.
\end{cor}

\ \\

We will prove this corollary at the end of Section 3.

\section{Computing Piecewise Step-Polynomials for Parametric Polytopes}

\label{s:explicit}

In this section, we recall the main elements
of the proof of Proposition~\ref{p:explicit} from \cite{BP99} and \cite{Verdoolaege2007parametric}.
That is,
given a parametric polytope $P\subset\Q^{n}\times\Q^{d}$, we define
$P_{\vec s}=\{\vec t\in\Q^{d} \mid (\vec s,\vec t)\in P\}$, for
$\vec s\in\Z^n$, and we want to compute
\[c(\vec s)=\#(P_{\vec s}\cap\Z^{d})\]
as a piecewise step-polynomial.
We demonstrate each step with a running example and
formulate an extended version
of the final step for use in Section~\ref{s:equivalence}.

\begin{example}
\label{ex:enumeration1}
Consider the parametric polytope
$$
P =
   \left\{
       (\vec s, \vec t) \in \Q^2 \times \Q^2 \ \ \bigg|\ \
       \left(
       \begin{matrix}
       -1 & 2 \\
       1 & -1 \\
       0 & 0 \\
       0 & 0
       \end{matrix}
       \right)
       \vec s +
       \left(
       \begin{matrix}
       1 & -2 \\
       -1 & 1 \\
       1 & 0 \\
       0 & 1
       \end{matrix}
       \right)
       \vec t \ge \vec 0
   \right\}
.
$$
We want to compute a piecewise step-polynomial representation of
$$
c(\vec s) = \#(P_{\vec s}\cap\Z^2)=\#\left\{\vec t\in\Z^{2} \mid
(\vec s,\vec t)\in P\right\} .
$$
\end{example}\boxend

Our main tool will be a slightly different sort of generating function than we have been using. If
$S\subset\Z^d$ is a set of integer vectors, then define its generating
function to be
\[
f(S;\x)=\sum_{\vec t\in S}\x^{\vec t}=\sum_{(t_1,\ldots,t_d)\in S}x_1^{t_1}x_2^{t_2}\cdots x_d^{t_d}.
\]
In our previous notation, this is the generating function for $c(\vec t)$ such that $c(\vec t)=1$ for $t\in S$ and $c(\vec t)=0$ otherwise.

Our proof of Proposition~\ref{p:explicit} will have two main steps.

\begin{itemize}
\item First, we will calculate the generating function
$f(P_{\vec s}\cap\Z^d;\x)$ as a rational generating function, and we
will examine how it changes as $\vec s$ varies (Propositions \ref{p:chambers} and \ref{p:core}).\\
\item Second, we will calculate
\[c(\vec s)=\#(P_{\vec s}\cap\Z^d)=f(P_{\vec s}\cap\Z^d;\vec 1),\]
by appropriately substituting $\x=\vec 1$.
\end{itemize}

In order to calculate the generating function $f(P_{\vec s}\cap\Z^d;\x)$, it is necessary to know what the vertices of $P_{\vec s}$ are.

\begin{figure}
\intercol=0.8cm
\begin{center}
\begin{minipage}{0cm}
\begin{xy}
<\intercol,0pt>:<0pt,\intercol>:: \POS(-1,0)\ar(10,0)
\POS(0,-1)\ar(0,5) \POS(10,0)*+!DR{s_1} \POS(0,5)*+!UR{s_2}
\POS(0,0)\ar@{-}@[|(2)](10,5) \POS(0,0)\ar@{-}@[|(2)](5,5)
\POS(3,5)\ar@{--}(3,2) \POS(3,2)\ar@{--}(10,2)
\POS(3,4)*{\bullet}="p1" \POS(3,3)*{\bullet}="p2"
\POS(3,2)*{\bullet}="p3" \POS(4,2)*{\bullet}="p4"
\POS(5,2)*{\bullet}="p5" \POS(1,4.5)*{Q_3} \POS(7.5,4.5)*{Q_1}
\POS(9,1)*{Q_2} \POS(-4,4)*[*0.6]\xybox{
<0.8\intercol,0pt>:<0pt,0.8\intercol>:: \POS(-5,0)\ar(5,0)
\POS(0,-1)\ar(0,5) \POS(5,0)*+!DR{t_1} \POS(0,5)*+!UR{t_2}
\POS(-5,0)\ar@{--}(5,5) \POS(-2,-1)\ar@{--}(4,5)
\POS@i@={(0,1),(0,2.5),(3,4),(0,1)},{0*[|(2)]\xypolyline{}}
\POS(0,1)*{\bullet}*+!DR{\vec v_5}
\POS(0,2.5)*{\bullet}*+!DR{\vec v_1}
\POS(3,4)*{\bullet}*+!DR{\vec v_6}
\POS(0,0)*{\circ}*+!DL{\vec v_2} \POS(-1,0)*{\circ}*+!DR{\vec v_3}
\POS(-5,0)*{\circ}*+!UL{\vec v_4} }*\frm{-}="s1"
\POS(-4,0)*[*0.6]\xybox{ <0.8\intercol,0pt>:<0pt,0.8\intercol>::
\POS(-5,0)\ar(5,0) \POS(0,-1)\ar(0,5) \POS(5,0)*+!DR{t_1}
\POS(0,5)*+!UR{t_2} \POS(-5,-1)\ar@{--}(5,4)
\POS(-1,-1)\ar@{--}(5,5)
\POS@i@={(0,0),(0,1.5),(3,3),(0,0)},{0*[|(2)]\xypolyline{}}
\POS(0,1.5)*{\bullet}*+!DR{\vec v_1}
\POS(3,3)*{\bullet}*+!DR{\vec v_6}
\POS(0,0)*{\bullet}*+!UL{\vec v_2 = \vec v_3 = \vec v_5}
\POS(-3,0)*{\circ}*+!UL{\vec v_4} }*\frm{-}="s2"
\POS(-4,-3)*[*0.6]\xybox{ <0.8\intercol,0pt>:<0pt,0.8\intercol>::
\POS(-5,0)\ar(5,0) \POS(0,-4)\ar(0,3) \POS(5,0)*+!DR{t_1}
\POS(0,3)*+!UR{t_2} \POS(-5,-2)\ar@{--}(5,3)
\POS(-3,-4)\ar@{--}(4,3)
\POS@i@={(0,0),(0,0.5),(3,2),(1,0),(0,0)},{0*[|(2)]\xypolyline{}}
\POS(0,-1)*{\circ}*+!UL{\vec v_5}
\POS(0,0.5)*{\bullet}*+!DR{\vec v_1}
\POS(3,2)*{\bullet}*+!DR{\vec v_6}
\POS(0,0)*{\bullet}*+!UR{\vec v_2}
\POS(1,0)*{\bullet}*+!UL{\vec v_3}
\POS(-1,0)*{\circ}*+!DR{\vec v_4} }*\frm{-}="s3"
\POS(2,-3)*[*0.6]\xybox{ <0.8\intercol,0pt>:<0pt,0.8\intercol>::
\POS(-5,0)\ar(5,0) \POS(0,-4)\ar(0,3) \POS(5,0)*+!DR{t_1}
\POS(0,3)*+!UR{t_2} \POS(-5,-2.5)\ar@{--}(5,2.5)
\POS(-2,-4)\ar@{--}(5,3)
\POS@i@={(0,0),(4,2),(2,0),(0,0)},{0*[|(2)]\xypolyline{}}
\POS(0,-2)*{\circ}*+!UL{\vec v_5}
\POS(0,0)*{\bullet}*+!DR{\vec v_1 = \vec v_2 = \vec v_4}
\POS(4,2)*{\bullet}*+!DR{\vec v_6}
\POS(2,0)*{\bullet}*+!UL{\vec v_3} }*\frm{-}="s4"
\POS(8,-3)*[*0.6]\xybox{ <0.8\intercol,0pt>:<0pt,0.8\intercol>::
\POS(-5,0)\ar(5,0) \POS(0,-4)\ar(0,3) \POS(5,0)*+!DR{t_1}
\POS(0,3)*+!UR{t_2} \POS(-5,-3)\ar@{--}(5,2)
\POS(-1,-4)\ar@{--}(5,2)
\POS@i@={(1,0),(5,2),(3,0),(0,0)},{0*[|(2)]\xypolyline{}}
\POS(0,-3)*{\circ}*+!UL{\vec v_5}
\POS(0,-0.5)*{\circ}*+!UL{\vec v_1}
\POS(5,2)*{\bullet}*+!DR{\vec v_6}
\POS(0,0)*{\circ}*+!DR{\vec v_2}
\POS(3,0)*{\bullet}*+!UL{\vec v_3}
\POS(1,0)*{\bullet}*+!DR{\vec v_4} }*\frm{-}="s5"
\POS"p1"\ar@{.>}@/_5ex/"s1" \POS"p2"\ar@{.>}@/_2ex/"s2"
\POS"p3"\ar@{.>}@/_2ex/"s3" \POS"p4"\ar@{.>}@/^2ex/"s4"
\POS"p5"\ar@{.>}@/^2ex/"s5"
\end{xy}
\end{minipage}
\end{center}
\caption{The decomposition and the vertices
of the parametric polytope from Example~\ref{ex:chambers}.}
\label{fig:chambers}
\end{figure}
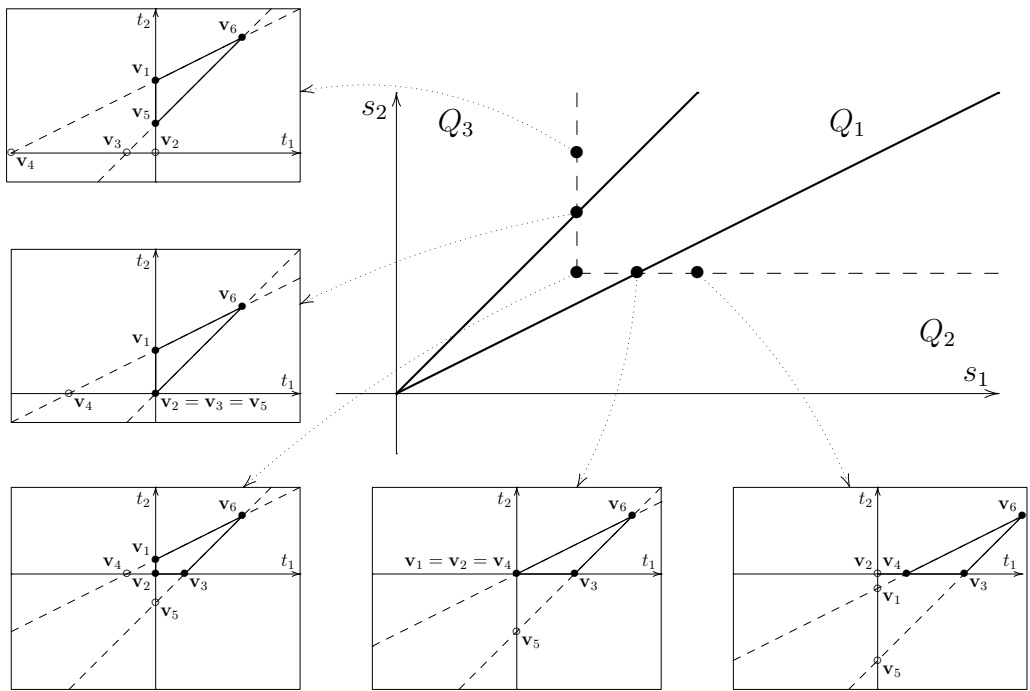

\begin{example}
\label{ex:chambers}  Consider the parametric polytope $P$ from Example~\ref{ex:enumeration1}.

For a given $\vec s$, the vertices of $P_{\vec s}$ can be obtained as the
intersections of pairs of facets of $P_{\vec s}$.  The
facets $t_1 = 0$ and $t_1 - 2 t_2 = s_1-2s_2$, for example,
intersect at the point $\vec v_1 = (0, -s_1/2 + s_2)$. This point is not always \emph{active}, that is, actually a vertex of $P_{\vec s}$. It is active exactly when $ 2s_2 \ge s_1 \ge 0$ (for
all other values of $\vec s$, $\vec v_1\notin P_{\vec s}$). We
similarly find the vertices $\vec v_2 = (0, 0)$,
$\vec v_3 = (s_1 - s_2, 0)$, $\vec v_4 = (s_1 - 2 s_2, 0)$,
$\vec v_5 = (0, -s_1 + s_2)$ and $\vec v_6 = (s_1, s_2)$, active
on the domains $2 s_2 \ge s_1 \ge s_2$, $s_1 \ge s_2 \ge 0$, $s_1
\ge 2 s_2 \ge 0$, $s_2 \ge s_1 \ge 0$, and $s_1, s_2 \ge 0$,
respectively. Combining all of the inequalities, we have the regions
\begin{eqnarray*}
Q_1 & = & \{\, \vec s \mid  2 s_2 \ge s_1 \ge s_2 \,\}
\\
Q_2 & = & \{\, \vec s \mid  s_1 \ge 2 s_2 \ge 0 \,\}
\\
Q_3 & = & \{\, \vec s \mid  s_2 \ge s_1 \ge 0 \,\} .
\end{eqnarray*}
For $\vec s\in Q_1$, the polyhedron $P_{\vec s}$ has active vertices $\vec v_1,\vec v_2, \vec v_3, \vec v_6$; for $\vec s\in Q_2$, it has active vertices $\vec v_3,\vec v_6,\vec v_4$; and for $\vec s\in Q_3$, it has active vertices $\vec v_1,\vec v_5,\vec v_6$. On the boundary of the $Q_i$, there is more than one possible description of the vertices (any is fine).

Figure~\ref{fig:chambers} shows the decomposition,
the vertices active in each $Q_i$, and the evolution
of the vertices as the value of $\vec s$ changes.
\end{example}
\boxend

As the example suggests, and as shown in \cite{Verdoolaege2007parametric}, we can find polyhedra $Q_i$ such that the $\interior(Q_i)$ partition $\Q^n$ and, for any $\vec s$ in the relative interior of a given polyhedron $Q_i$, the polytopes $P_{\vec s}$ will have a fixed set of
vertices given by affine transformations of $\vec s$ (where an
affine transformation $T:\Q^n\rightarrow \Q^d$ is given by $T(\vec
s)=T'(\vec s)+\vec v$ such that $T'$ is a linear transformation and
$\vec v\in\Q^d$). These $Q_i$ will be the pieces of our
piecewise step-polynomial.  This is the content of the following proposition.

\begin{prop}[Decomposition]
\label{p:chambers} Fix $d$ and $n$. There exists a polynomial time
algorithm, which, given a parametric polytope $P \subset \Q^{n}
\times \Q^{d}$, finds polyhedra $Q_i$ whose relative interiors partition $\Q^n$, and, for
each $i$, computes a collection of affine transformations $T_{i1},
T_{i2},\ldots, T_{im_i}:\Q^{n}\rightarrow\Q^{d}$, such that, for
$\vec s\in \interior{Q_i}$, the vertices of $P_{\vec s}$ are $T_{i1}(\vec
s),T_{i2}(\vec s),\ldots,T_{im_i}(\vec s)$.
\end{prop}

Algorithms to compute the parametric vertices and the chambers can
be found in~\cite{Loechner97parameterized} and~\cite{CL98}
respectively. A proof of the polynomial time complexity is given
in~\cite{Verdoolaege2007parametric}.

Now we can concentrate on computing $f(P_{\vec s}\cap\Z^d;\x)$, given that $\vec s$ is in the relative interior of a particular $Q_i$. As a first step, we examine how to compute the generating function
of an easy set: the integer points in a \emph{unimodular cone}.
The general case of a polyhedron is based on a reduction to these
unimodular cones.
\begin{definition}
Let $\vec c_1,\vec c_2,\ldots,\vec c_d\in\Z^d$ be a basis for the
lattice $\Z^d$, and let $\beta_i\in \Q$, for $1\le i\le d$.  We
define the \emph{rational unimodular cone}
\[K=\left\{\,\x\in\Q^d \mid  \langle \vec c_i,\vec x\rangle\le \beta_i
\text{ for }1\le i\le d\,\right\}.\]
\end{definition}
This cone may have a vertex which is not at the origin.  Let
$\vec u_1,\vec u_2,\ldots,\vec u_d$ be the \emph{negative dual basis} of
$\Z^d$, so that
\[\langle \vec u_i,\vec c_j\rangle=\left\{%
\begin{array}{ll}
    -1, & \hbox{if }i=j \\
    0, & \hbox{if }i\ne j
.
\end{array}%
\right.\]
If $\beta_i$ were zero, for all $i$, then $K$ would be
the cone with vertex at the origin defined by
\[
K=\big\{\,\lambda_1\vec u_1+\lambda_2\vec u_2+\cdots+\lambda_d\vec u_d
\ \mid \
 \vec\lambda\ge \vec0\,\big\},
\]
we would have that
\[
K\cap \Z^d
=\big\{\,\lambda_1\vec u_1+\lambda_2\vec u_2+\cdots+\lambda_d\vec u_d
\ \mid \
 \vec\lambda \in \Z_{\ge 0}^d\,\big\},
\]
and therefore
\[f(K\cap\Z^d;\x)=
  \frac{1}{(1-\x^{\vec u_1})(1-\x^{\vec u_2})\cdots(1-\x^{\vec u_d})}.\]
In the general case, where the $\beta_i$ are not necessarily zero,
we have that
\begin{equation}
\label{eq:p}
f(K\cap\Z^d;\x)=
  \frac{\x^{\vec p}}{(1-\x^{\vec u_1})(1-\x^{\vec u_2})\cdots(1-\x^{\vec u_d})},
\end{equation}
where $\vec p=-\sum_{i=1}^d\floor{\beta_i}\vec u_i$ (see
\cite{BP99}).  This greatest integer function in the definition of
$\vec p$ is where the greatest integer function in our
step-polynomial will come from. Note also that the denominator of
this generating function does not depend on the $\beta_i$, only on
the $\vec c_i$.

We want to reduce our problem, which is finding the generating
function $f(P_{\vec s}\cap\Z^d;\x)$ where $P_{\vec s}$ is a
polyhedron, to the easy problem of finding the generating function
for a unimodular cone.  We can first reduce to the case of (not necessarily unimodular) cones using Brion's Theorem
\cite{Brion88}, which states that the generating function of a polytope
is equal to the sum of the generating functions of its \emph{vertex cones}.
These vertex cones are formed by the supporting hyperplanes of the
polytope that intersect in a given vertex (see Figure~\ref{fig:C3} for
an example).
Next, we use Barvinok's unimodular decomposition \cite{Barvinok94} to write the generating function of each vertex cone
as a (signed) sum of generating functions of unimodular cones.

\begin{example}
\label{ex:generating}
Consider once more the parametric polytope $P$ from Examples~\ref{ex:enumeration1} and~\ref{ex:chambers}.
We want to compute the generating function of this parametric polytope.
Consider specifically region $Q_3$ from Example~\ref{ex:chambers}
with active vertices $\vec v_1 = (0, -s_1/2
+ s_2)$, $\vec v_5 = (0, -s_1 + s_2)$ and $\vec v_6 = (s_1, s_2)$.
The polytope corresponding to $\vec s = (3,4) \in C_3$ is shown in
Figure~\ref{fig:C3} together with the vertex cones, $\cone(P_{\vec s},\vec v_i)$, at each
active vertex. Brion's theorem tells us that
\begin{align*} f(P_{\vec s}\cap\Z^2;\vec x)=f\big(\cone&(P_{\vec s},\vec v_1)\cap\Z^2;\vec x\big)+f\big(\cone(P_{\vec s},\vec v_5)\cap\Z^2;\vec x\big)\\
&+f\big(\cone(P_{\vec s},\vec v_6)\cap\Z^2;\vec x\big).
\end{align*}
The vertex cones at $\vec v_5$ and
$\vec v_6$ are unimodular, but the one at $\vec v_1$ is not. We
therefore need to apply Barvinok's unimodular decomposition to $\cone(P_{\vec s},\vec v_1)$, which yields

\begin{equation}\label{eq:beforesubs} f\big(\cone(P_{\vec s},\vec v_1)\cap\Z^2;\vec x\big) =
\frac {\x^{(-2 \lfloor \frac {s_1} 2 - s_2 \rfloor + s_1 - 2 s_2,
            - \lfloor \frac {s_1} 2 - s_2 \rfloor )}}
      {(1 - \x^{(1,0)})(1 - \x^{(2,1)})}
-
\frac {\x^{(0, -\lfloor \frac {s_1} 2 - s_2 \rfloor)}}
      {(1 - \x^{(1,0)})(1 - \x^{(0,1)})}
.
\end{equation}
We refer to \cite{BP99,DHT03,Koeppe2006primal,Koeppe2007parametric}
for details on how to perform Barvinok's decomposition.
Table~\ref{tab:generating} lists the generating functions of all vertex cones.

\begin{figure}
\intercol=0.8cm
\begin{xy}
<\intercol,0pt>:<0pt,\intercol>:: \POS(-1,0)\ar(5,0)
\POS(0,-1)\ar(0,5) \POS(5,0)*+!DR{t_1} \POS(0,5)*+!UR{t_2}
\POS@i@={(0,1),(0,2.5),(3,4),(0,1)},{0*[|(2)]\xypolyline{}}
\POS(0,1)*{\bullet}*+!DR{\vec v_5}
\POS(0,2.5)*{\bullet}*+!DR{\vec v_1}
\POS(3,4)*{\bullet}*+!DR{\vec v_6} \POS(-5,2)*\xybox{
<0.7\intercol,0pt>:<0pt,0.7\intercol>::
\POS@i@={(0,-1),(0,2.5),(4,4.5),(4,-1)},{0*[grey]\xypolyline{*}}
\POS@i@={(0,-1),(0,2.5),(4,4.5)},{0*[|(2)]\xypolyline{}}
\POS(-1,0)\ar(4.5,0) \POS(0,-1)\ar(0,5)
\POS(0,2.5)*{\bullet}*+!DR{\vec v_1} } \POS(6.5,4)*\xybox{
<0.7\intercol,0pt>:<0pt,0.7\intercol>::
\POS@i@={(-1,2),(3,4),(-1,0)},{0*[grey]\xypolyline{*}}
\POS@i@={(-1,2),(3,4),(-1,0)},{0*[|(2)]\xypolyline{}}
\POS(-1,0)\ar(4.5,0) \POS(0,-1)\ar(0,5)
\POS(3,4)*{\bullet}*+!DR{\vec v_6} } \POS(6.5,-1)*\xybox{
<0.7\intercol,0pt>:<0pt,0.7\intercol>::
\POS@i@={(0,4.5),(0,1),(3.5,4.5)},{0*[grey]\xypolyline{*}}
\POS@i@={(0,4.5),(0,1),(3.5,4.5)},{0*[|(2)]\xypolyline{}}
\POS(-1,0)\ar(4.5,0) \POS(0,-1)\ar(0,5)
\POS(0,1)*{\bullet}*+!DR{\vec v_5} }
\end{xy}
\caption{$P_{(3,4)}$ and its vertex cones}
\label{fig:C3}
\end{figure}
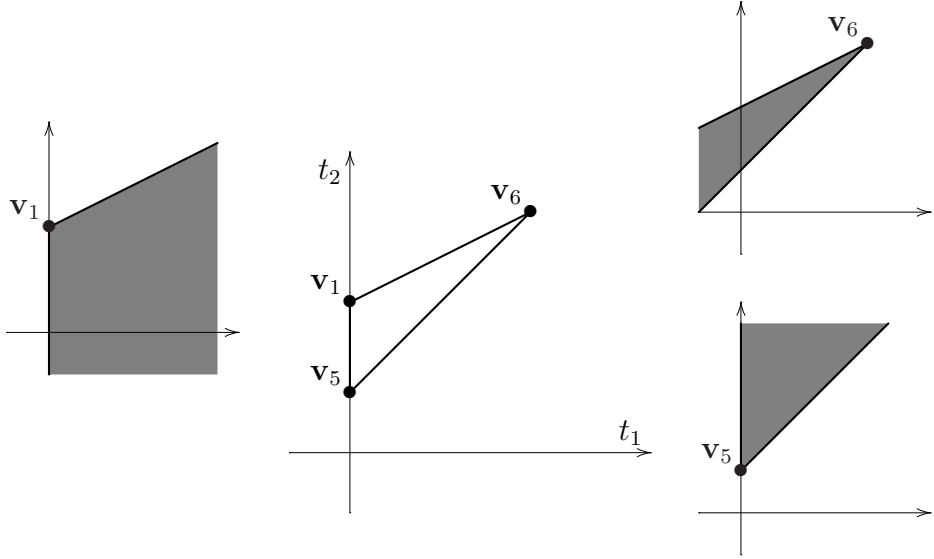

\begin{table}
$$
\begin{array}{|c|c|}
\hline \hbox{Vertex $\vec v_i$} &
\hbox{$f\big(\cone(P_{\vec s},\vec v_i)\cap\Z^d;\vec x\big)$} \\
\hline \hline \vec v_1 = (0, -s_1/2 + s_2) &
\frac {\x^{(-2 \lfloor \frac {s_1} 2 - s_2 \rfloor + s_1 - 2 s_2,
            - \lfloor \frac {s_1} 2 - s_2 \rfloor )}}
      {(1 - \x^{(1,0)})(1 - \x^{(2,1)})}
-
\frac {\x^{(0, -\lfloor \frac {s_1} 2 - s_2 \rfloor)}}
      {(1 - \x^{(1,0)})(1 - \x^{(0,1)})}
\\
\hline \vec v_2 = (0, 0) &
\frac {\x^{(0, 0)}}
      {(1 - \x^{(0,1)})(1 - \x^{(1,0)})}
\\
\hline \vec v_3 = (s_1 - s_2, 0) &
\frac {\x^{(s_1 - s_2, 0)}}
      {(1 - \x^{(1,0)})(1 - \x^{(-1,0)})}
\\
\hline \vec v_4 = (s_1 - 2 s_2, 0) &
\frac {\x^{(s_1 - 2 s_2, 0)}}
      {(1 - \x^{(2,1)})(1 - \x^{(1,0)})}
\\
\hline \vec v_5 = (0, -s_1 + s_2) &
\frac {\x^{(0, -s_1 + s_2)}}
      {(1 - \x^{(0,1)})(1 - \x^{(1,1)})}
\\
\hline \vec v_6 = (s_1, s_2) &
\frac {\x^{(s_1, s_2)}}
      {(1 - \x^{(-2,-1)})(1 - \x^{(-1,-1)})}
\\
\hline
\end{array}
$$
\caption{The generating function of each vertex cone}
\label{tab:generating}
\end{table}

\end{example}
\boxend

When we do this in general, the end result is the following proposition, a rephrasing of
Theorem 4.4 of \cite{BP99}.

\begin{prop}
\label{p:core} Fix $d$.  There exists a polynomial time algorithm,
which, given a parametric polyhedron $P \subset \Q^n\times\Q^d$ and
a polyhedral region $Q$ such that for $\vec s\in Q$ the vertices of
$P_{\vec s}=\{\vec t\in\Q^d\ \mid\ (\vec s, \vec t)\in P\}$ are
given by affine transformations $T_1(\vec s),T_2(\vec
s),\cdots,T_m(\vec s)$, computes the generating function
$$
f(P_{\vec s} \cap \Z^d; \x)=\sum_{i\in I}\varepsilon_i
\frac{\x^{\vec p_i({\vec s})}}
     {(1-\x^{\vec b_{i1}})(1-\x^{\vec b_{i2}})\cdots(1-\x^{\vec b_{id}})},
$$
where $\varepsilon\in \{\, -1, 1 \,\}$,
$\vec b_{ij}\in\Z^d\setminus\{\vec 0\}$, and each coordinate of
$\vec p_i(\vec s):\Z^n\rightarrow\Z^d$ is a step-polynomial of
degree one, for each $i$.
\end{prop}

Now that we know how to compute $f(P_{\vec s}\cap\Z^d; \x)$, all that
remains is to evaluate it at $\vec x = \vec 1$.

\begin{example}
\label{ex:enumeration}
Consider once more the parametric polytope $P$ from Examples~\ref{ex:enumeration1}, \ref{ex:chambers},
and \ref{ex:generating}.

We have already computed $f\big(P_{\vec s}\cap \Z^2;
\vec x)$, and
we now compute the value $f\big(P_{\vec s} \cap \Z^2;
\vec 1)$.  We cannot simply plug in $\vec x = \vec 1$, because $\vec 1$ is a pole of some of the rational functions.  Instead, we make a suitable substitution, in this case $\x = (t+1, t+1)$ (chosen carefully so that none of the denominators become identically zero), and take the limit as $t$ approaches zero. To compute this limit, we can simply compute, for each term in the sum constituting $f\big(P_{\vec s}\cap \Z^2;
t+1,t+1)$, the constant term in the Laurent series expansion at $t=0$.

For example, substituting $\vec x=(t+1,t+1)$ into the second term in
(\ref{eq:beforesubs}), we obtain
$$
- \frac {(1+t)^{-\lfloor \frac {s_1} 2 \rfloor + s_2}}
        {(1- (1+t)) (1-(1+t))}
.
$$
Since the denominator, in this case, is exactly $t^2$, the constant
term in the Laurent expansion is simply the coefficient of $t^2$ in
the expansion of the numerator, i.e.,
$$
- \frac { (-\lfloor \frac {s_1} 2 \rfloor + s_2)
          (-\lfloor \frac {s_1} 2 \rfloor + s_2 -1 ) } 2
=
- \frac 1 2
\left\lfloor \frac {s_1} 2 \right\rfloor^2
- \frac {s_2^2} 2
+
\left\lfloor \frac {s_1} 2 \right\rfloor s_2
- \frac 1 2
\left\lfloor \frac {s_1} 2 \right\rfloor
+
\frac {s_2} 2
.
$$
The other terms are handled similarly. Note that this is the place where the step polynomials show up in full force. The contribution of each vertex cone to the constant term of the Laurent expansion is
listed in Table~\ref{tab:constant}. The final step-polynomial in
each chamber is computed using Brion's Theorem as the sum of the
appropriate step-polynomials from this table. The final result is
shown in Figure~\ref{fig:result}.

\begin{table}
$$
\begin{array}{|c|c|}
\hline \hbox{Vertex $\vec v_i$} &
\hbox{$f\big(\cone(P_{\vec s},\vec v_i)\cap\Z^d;\vec 1\big)$ (if $v_i$ is active)} \\
\hline \hline \vec v_1 = (0, -s_1/2 + s_2) & \frac {s_1^2} 6 +
\frac {s_1 s_2} 3 - s_1 \left\lfloor \frac {s_1} 2 \right\rfloor -
\frac {s_1} 2 - \frac {s_2} 3 + \left\lfloor \frac {s_1} 2
\right\rfloor^2 + \left\lfloor \frac {s_1} 2 \right\rfloor + \frac 2
9
\\
\cline{2-2}
&
- \frac 1 2
\left\lfloor \frac {s_1} 2 \right\rfloor^2
- \frac {s_2^2} 2
+
\left\lfloor \frac {s_1} 2 \right\rfloor s_2
- \frac 1 2
\left\lfloor \frac {s_1} 2 \right\rfloor
+
\frac {s_2} 2
\\
\hline \vec v_2 = (0, 0) & 0
\\
\hline \vec v_3 = (s_1 - s_2, 0) & - \frac {s_1^2} 4 + \frac {s_1
s_2} 2 - \frac {s_2^2} 4 + \frac 1 8
\\
\hline \vec v_4 = (s_1 - 2 s_2, 0) & \frac {s_1^2} 6 - \frac {2 s_1
s_2} 3 - \frac {s_1} 2 + \frac {2 s_2^2} 3 + s_2 + \frac 2 9
\\
\hline \vec v_5 = (0, -s_1 + s_2) & \frac {s_1^2} 4 - \frac {s_1
s_2} 2 + \frac {s_1} 2 + \frac {s_2^2} 4 - \frac {s_2} 2 + \frac 1 8
\\
\hline \vec v_6 = (s_1, s_2) & \frac {s_1^2} {12} + \frac {s_1 s_2} 6
+ \frac {s_1} 2 + \frac {s_2^2}{12} + \frac {s_2} 2 + \frac {47}{72}
\\
\hline
\end{array}
$$
\caption{The contribution of each vertex cone to the constant
term of the Laurent expansion}
\label{tab:constant}
\end{table}

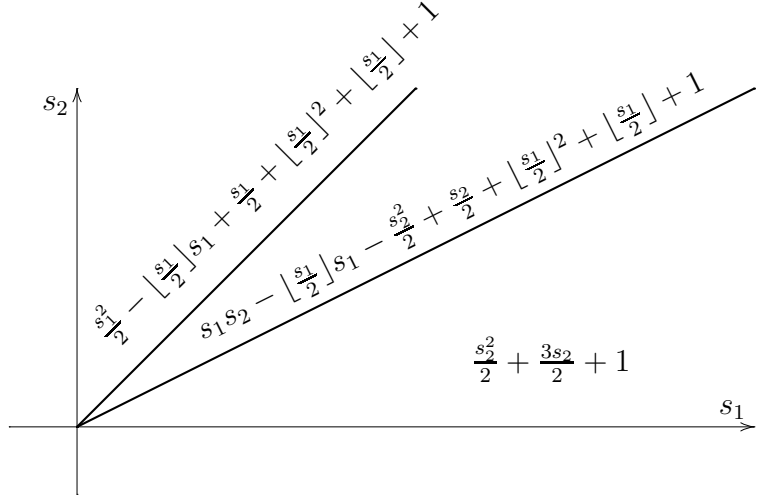
\begin{figure}
\intercol=0.9cm
\begin{center}
\begin{minipage}{0cm}
\begin{xy}
<\intercol,0pt>:<0pt,\intercol>::
\POS(-1,0)\ar(10,0)
\POS(0,-1)\ar(0,5)
\POS(10,0)*+!DR{s_1}
\POS(0,5)*+!UR{s_2}
\POS(0,0)\ar@{-}@[|(2)](10,5)
\POS(0,0)\ar@{-}@[|(2)](5,5)
\POS(2.8,3.8)*[@v(1,1)]{\frac{s_1^2}2 - \lfloor \frac{s_1}2 \rfloor s_1
+ \frac {s_1} 2
+ \lfloor \frac{s_1}2 \rfloor^2
+ \lfloor \frac{s_1}2 \rfloor +1}
\POS(5.5,3.3)*[@v(2,1)]{s_1 s_2 - \lfloor \frac{s_1}2 \rfloor s_1
- \frac {s_2^2}2 + \frac {s_2} 2
+ \lfloor \frac{s_1}2 \rfloor^2
+ \lfloor \frac{s_1}2 \rfloor +1}
\POS(7,1)*{\frac {s_2^2}2 + \frac {3s_2}2 +1}
\end{xy}
\end{minipage}
\end{center}
\caption{The enumerator of $P_{\vec s}$, a step-polynomial
in each chamber}
\label{fig:result}
\end{figure}

\end{example}
\boxend

In general, we use the following
lemma, which is more general than strictly needed here, but which
allows for an incremental computation as discussed after the lemma
and which we will also need in Section~\ref{s:equivalence}.
The lemma is a special case %pun removed :(
of the monomial substitution theorem \cite[Theorem 2.6]{BW03}. We
provide a slightly different proof, which lends itself more easily
to an implementation. It is an extension of an idea from
\cite{DHT03}, which is in itself a variation of the idea used in
\cite{Barvinok94}.

\begin{lemma}[Specialization]
\label{l:specialization} Let us fix $k$. There exists a polynomial
time algorithm which, given a rational generating function $f(\x)$
of the form (\ref{GFForm}) and an $m$ with $0 \le m \le d$ such that
$g(\vec z) := f(z_1, \ldots, z_m, 1, \ldots, 1)$ is an analytic
function on some nonempty open subset of $\C^m$, computes $g(\vec
\z)$ in the same form, i.e.,
\begin{equation}
\label{eq:specialized} g(\vec z) = \sum_{i' \in I'} \beta_{i'}
\frac{\vec z^{\vec q_{i'}}}
     {(1-\vec z^{\vec d_{i'1}})(1-\vec z^{\vec d_{i'2}})\cdots(1-\vec z^{\vec d_{i'k_{i'}}})},
\end{equation}
where $k_{i'} \le k$, $\z\in \C^m$, $\beta_{i'}\in\Q$,
$\vec q_{i'}\in\Z^m$, and $\vec d_{i'j'}\in\Z^m\setminus \{\vec 0\}$.

Furthermore, if the vectors $\vec b_{ij}$ and the numbers
$\alpha_i$ in (\ref{GFForm}) are fixed, but the vectors $\vec p_i$
vary, then the vectors $\vec d_{i'j'}$ are fixed, $\vec q_{i'}$
each differ by a constant vector from some $\vec p_i$, and
$\beta_{i'}$ are each a polynomial of degree at most $k$ in the
coordinates of some $\vec p_i$.
\end{lemma}

\begin{proof}
The case $m = d$ is trivial, so we will assume $m < d$. Note that we
cannot simply plug in the values $1$, since $(z_1, \ldots, z_m, 1,
\ldots, 1)$ may be a pole of some of the terms in (\ref{GFForm}). In
fact, if $m = 0$, then it will be a pole of all those terms. We must
take an appropriate limit as $(x_{m+1},\ldots,x_n)$ approaches
$(1,\ldots,1)$.  Consider
$$
h(t) = f(z_1, \ldots, z_m, (1+t)^{\lambda_1}, \ldots,
(1+t)^{\lambda_{d-m}}) ,
$$
as a function of $t$ only, where $\vec \lambda\in\Z^{d-m}$ is such
that for each $i \in I$ either $(b_{i1}, \ldots, b_{im}) \ne \vec 0$
or $\langle (b_{i,m+1}, \ldots, b_{id}), \vec \lambda \rangle \ne 0
$. Such a $\vec \lambda$ can be found in polynomial time by choosing
an appropriate point from the ``moment curve'' as in \cite[Algorithm
5.2]{BP99}. Then $g(\vec z)$ is simply the constant term in the
Laurent power series expansion of $h(t)$ about $t=0$.  This is the
sum of the constant terms in the Laurent power series expansions of
$$
h_i(t) = \alpha_i \frac{ \vec z^{\vec p_i'} (t+1)^{q_i}} {(1-\vec
z^{\vec b'_{i1}}(t+1)^{v_{i1}})
 (1-\vec z^{\vec b'_{i2}}(t+1)^{v_{i2}})\cdots
 (1-\vec z^{\vec b'_{ik_i}}(t+1)^{v_{ik_i}})},
$$
where, for $\vec v \in \{\, \vec p_i, \vec b_{ij} \,\}$, we write
$\vec v'$ for the first $m$ components of $\vec v$ and $\vec v''$
for the remaining $d-m$ components, and we let $q_i=\langle
\vec p'',\vec 1\rangle$ and $v_{ij}=\langle \vec b''_{ij}, \vec
1\rangle$.

Consider a particular $h_i(t)$. Let $r$ be the number of factors with
$v_{ij} \ne 0$ but $\vec b'_{ij} = \vec 0$.  Then $h_i(t)$ has a
pole of order $r$ at $t=0$.  Therefore, we must compute the
coefficient of $t^r$ in the Taylor series expansion of $t^rh_i(t)$,
which is analytic at $t=0$.
%Partition the factors in the denominator based on the variables that
%appear in each factor. Let $q$ be the number of factors such that
%$v_{il} = 0$, $r$ the number of factors with $v_{il} \ne 0$ but
%$\vec b'_{il} = \vec 0$ and $s = k - q - r$ the number of factors
%with both $v_{il} \ne 0$ and $\vec b'_{il} \ne \vec 0$. I.e.,
%$$
%h_i(t) = C(\vec z) \frac{(t+1)^{q_i}} { \prod_r ((t+1)^{\beta_j} -
%1) \prod_s ((t+1)^{\gamma_j} - \vec z^{\vec c_j}) } = C(\vec z)
%\frac {P(t)}{t^r \prod Q_j(t; \vec z)} ,
%$$
%where $C(\vec z)$ is a symbolic constant that collects the factors
%in both numerator and denominator that are independent of $t$ and
%where we have multiplied both numerator and denominator with either
%$-1$ or $-\vec z^{\vec c_j} = -\vec z^{\vec b'_{ij'}}$ for some
%$j'$. In the special that $r = 0$ (in particular, this requires $m
%\ne 0$), then $t=0$ is not a pole and we can simply plug in $0$ for
%$t$. We obtain a single term in the sum of the form
%(\ref{eq:specialized}) with $k' = k$.

%In general we need to be more careful and compute the coefficient of
%$t^r$ in $P(t)/ \prod Q_j(t; \vec z)$.

Following \cite{DHT03} we use the technique outlined in
\cite[241--247]{Henrici1974} (where it is applied to compute the
residue of a function, i.e., the coefficient of the term $t^{-1})$.
Let $t^rh_i(t)=\frac{P(t)}{Q(t)}$, where $P$ and $Q$ are
polynomials. To compute the coefficients $c_j$ in
$$
\frac{P(t)}{Q(t)} =: c_0 + c_1 t + c_2 t^2 + \cdots
,
$$
expand $P(t)$ and $Q(t)$ as
\begin{eqnarray*}
P(t) & =: & a_0 + a_1 t + a_2 t^2 + \cdots \\
Q(t) & =: & b_0 + b_1 t + b_2 t^2 + \cdots
\end{eqnarray*}
and apply the recurrence relation
$$
c_j = \frac 1 {b_0} \left( a_j - \sum_{i=1}^j b_i c_{j-i} \right)
.
$$
Note that we only need to keep track of the first $r+1$ coefficients
of $P(t)$ and $Q(t)$, and so this may be done in polynomial time.
Examining the recursive process, we see that the lemma follows.
%First divide by the first $r$ factors, which are independent of
%$\vec z$. The constant terms of the remaining factors are of the
%form $b_0 = 1 - \vec z^{\vec c_j}$. Only expressions of this kind
%will ever appear in a denominator. After the first division, the
%largest power in a denominator of $c_j^{[1]}$ is $b_0^{r+1}$. Each
%subsequent division increases the total power of all factors in the
%denominator by one. The total power of factors in the denominator of
%$c_j^{[s]}$ will therefore be $r+s$ and so remains constant. The
%number of terms in $c_j^{[s]}$ is also clearly bounded by a constant
%$s(k)$, i.e., $|I'| \le s(k) |I|$.
%
%Note that based on the binomial theorem
%$$
%(1+t)^r = \sum_k {r \choose k} t^k
%$$
%the coefficients of the numerator $P(t)$ are polynomial expressions
%in $q_i$, which is itself a linear combination of the coefficients
%of $\vec p_i$. The maximal degree of such a polynomial expression
%is $r \le k$. The lemma is proved.
\end{proof}

\noindent{\bf Remark on the implementation of Lemma
\ref{l:specialization}:} Note that as argued by \cite{DHT03}, a
$\vec \lambda$ from the moment curve may not be the most appropriate
choice to use in an implementation since it is likely to have large
coefficients. They therefore propose to construct a random vector
with small coefficients and check whether $\langle \vec b''_{ij},
\vec \lambda\rangle \ne 0$ for all $i$ and $j$. (Or rather $\langle
\vec b_{ij}, \vec \lambda\rangle \ne 0$, since $m=0$ in their
case.) Only after a fixed number of failed attempts would the
implementation fall back onto the moment curve.

Both of these strategies have the disadvantage however that
all the terms in (\ref{GFForm}) need to be available before the
constant term of the first term can be computed.
This may induce a large memory bottleneck.
The authors of \cite{DHT03} have therefore also implemented an alternative
strategy where a random vector with larger coefficients is constructed
at the beginning of the computation.  If the coefficients are large
enough, then the probability of having constructed an incorrect
vector is close to zero.
The disadvantage of this technique is that the coefficients are
larger and that the computation has to be redone completely
in the unlikely event the vector was incorrect.

We propose a different strategy which does not require all terms to
be available, nor does it require the use of large coefficients. We
simply repeatedly apply Lemma~\ref{l:specialization} for $m'$ from
$d-1$ down to $m$. In each application, we can simply use $\lambda =
1$, which is known to be valid in any case. \boxend

%\begin{example}
%Consider the rational generating function
%$$
%f(z,x) =
%\frac {x^4 z^{-1}}{(1-x)(1-x z^{-1})}
%+
%\frac {z^2}{(1-z^{-1})(1-x z^{-1})}
%-
%\frac {z^{-1}}{(1-x)(1-z^{-1})}
%.
%$$
%We want to compute $g(z) = f(z,1)$.
%Substituting $x = t + 1$, we obtain
%$$
%f(z,t+1) =
%\frac {(t+1)^4}{t(t + (1 - z))}
%+
%\frac {z^2}{(1-z^{-1})(1-(t+1) z^{-1})}
%+
%\frac {z^{-1}}{t(1-z^{-1})}
%.
%$$
%In the notation of the proof, the final term has $r=1$
%and $s=0$.  Since the coefficient of $t^1$ in the numerator
%is $0$, the contribution of this term is $0$.
%The second term has $r = 0$ and so we can simply plug in $1$
%to obtain
%$$
%\frac {z^2}{(1-z^{-1})^2}
%.
%$$
%The first term has $r=1$ and $s=1$.
%The numerator is
%$$
%(t+1)^4 \equiv 1 + 4 t \mod {t^2}
%,
%$$
%i.e., $a_0 = 1$ and $a_1 = 4$,
%while the denominator has $b_0 = 1-z$ and $b_1$.
%We find
%$$
%c_0 = \frac 1 {1 - z}
%\qquad\hbox{and}\qquad
%c_1 =  \frac 1 {1 - z}
%    \left(
%        4 - \frac 1 {1 - z}
%    \right)
%.
%$$
%Summing the three contributions, we obtain
%$$
%g(z) =
%0
%+
%\frac {z^2}{(1-z^{-1})^2}
%+
%\left(
%-
%\frac { 4 z^{-1}}{ 1 - z^{-1}}
%-
%\frac {  z^{-2}}{ (1 - z^{-1})^{-2}}
%\right)
%,
%$$
%where we applied (\ref{eq:identity}) to make all terms
%converge on a common neighborhood.
%\todo{Put (\ref{eq:identity}) before this use}
%\end{example}

We summarize the proof of Proposition~\ref{p:explicit}.

\begin{proof}[Proof of Proposition~\ref{p:explicit}]
\label{p:explicit:proof}
Given a parametric polytope $P\subset\Q^n\times\Q^d$, apply
Proposition~\ref{p:chambers} to obtain the decomposition
$\{\, Q_i \,\}$. For each region $Q_i$, apply
Proposition~\ref{p:core} to obtain the corresponding generating
function of $P_{\vec s}$, for $s\in Q_i$. The result is a collection
of polyhedral regions $Q_i$ such that, for $\vec s\in \interior(Q_i) \cap \Z^n$,
\[f(P_{\vec s}\cap\Z^d;\x)=\sum_j
\frac{\x^{\vec p_j(\vec
s)}}{(1-\x^{\vec u_{j1}})(1-\x^{\vec u_{j2}})\cdots
(1-\x^{\vec u_{jd}})},\] where $\vec u_{jl}\in\Z^d\setminus\{\vec 0\}$
and the coordinates of $\vec p_j:\Z^n\rightarrow\Z^d$ are piecewise
step-polynomials of degree one. All that remains is to use
Lemma~\ref{l:specialization} with $m=0$ to compute $c(\vec
s):=f(P_{\vec s}\cap\Z^d;\vec 1)$ as a step-polynomial in $\vec s$,
valid for $\vec s \in \interior{Q_i}$.
\end{proof}

\section{Equivalence of Rational Generating Functions
and Piecewise Step-Polynomials}

\label{s:equivalence}

In this section, we prove Theorem~\ref{ThmEquivalence}, that we may
convert between rational generating function and piecewise
step-polynomial representations in polynomial time. In both
directions, we reduce the problem to a set of counting problems to
which we apply either Proposition~\ref{p:explicit} or
Proposition~\ref{p:generating}. We first prove a special case of the
first half of Theorem \ref{ThmEquivalence}, as a corollary of
Proposition~\ref{p:explicit}.

\begin{cor}
\label{CorVPF} Fix $d$.  There is a polynomial time algorithm which,
given $\alpha\in\Q$, $\vec p\in\Z^n$,
$\vec a_1,\vec a_2,\ldots,\vec a_d\in\Z^n\setminus\{\vec 0\}$ and
given $\vec l\in\Z^n$ such that $\langle \vec l,\vec a_i\rangle\ne
0$ for all $i$, computes the piecewise step-polynomial
$c:\Z^n\rightarrow \Q$ such that
\[f(\x)=\alpha\frac{\x^{\vec p}}{(1-\x^{\vec a_1})(1-\x^{\vec a_2})\cdots(1-\x^{\vec a_d})}=
\sum_{\vec s\in\Z^n}c(\vec s)\x^{\vec s}\] is convergent on a
neighborhood of $\mathbf{e}^{\vec l}$.
\end{cor}

\begin{proof}
We may assume, without loss of generality, that $\langle \vec
l,\vec a_i\rangle < 0$ for all $i$. Otherwise, if $\langle \vec
l,\vec a_i\rangle
>0$ for some $i$, we would apply the identity
\begin{equation}
\label{eq:identity}
\frac{1}{1-\x^{\vec a_i}}=\frac{-\x^{-\vec a_i}}{1-\x^{-\vec a_i}}.
\end{equation}
It suffices to prove this corollary for $\alpha=1$ and  $\vec p=\vec
0$, because if $c'(\vec s)$ is a piecewise step-polynomial
representation of the generating function $g(\x)$, then $\alpha\cdot
c'(\vec s-\vec p)$ is a piecewise step-polynomial representation of
$\alpha\x^{\vec p}g(\x)$).
Note that $\alpha=1$ and $\vec p=\vec 0$
mean that $c(\vec s)$ is the vector partition function defined in Example
\ref{ex:VPFDef}.

We expand $f(\x)$ as a product of infinite geometric series,
\[f(\x)=\prod_{i=1}^{d}(1+\x^{\vec a_i}+\x^{2\vec a_i}+\cdots).\]
Then \[f(\mathbf{e}^{\vec l})=\prod_{i=1}^{d}(1+e^{\langle \vec
l,\vec a_i\rangle}+e^{2\langle \vec l,\vec a_i\rangle}+\cdots),\]
and this expansion is convergent on a neighborhood of
$\mathbf{e}^{\vec l}$, since $\langle \vec l,\vec a_i\rangle <0$.
We see that we are looking to compute the function
\[c(\vec s)=\#\big\{
    \vec\lambda=(\lambda_1,\lambda_2,\ldots,\lambda_d)\in\Z^d_{\ge
0} \ \mid \ \vec s=\lambda_1 \vec a_1+\lambda_2
\vec a_2+\cdots+\lambda_d \vec a_d\big\} .\]

Let $P$ be the parametric polytope
\[P=\big\{(\vec s, \vec \lambda)\in\Q^n \times \Q^d \ \mid \
\vec\lambda\ge\vec 0\text{ and }\vec s=\lambda_1
\vec a_1+\cdots+\lambda_d \vec a_d\big\} .\] Then
\[c(\vec s)=\#\left\{\vec\lambda\in\Z^d \mid (\vec s,\vec\lambda)\in P\right\},
\]
which can be computed as a piecewise step-polynomial using
Proposition~\ref{p:explicit}. The proof follows.
\end{proof}

\begin{example}
Consider the function
$$
f(\x) =
\frac 1
      {(1-\x^{(1,1)})(1-\x^{(2,1)})(1-\x^{(1,0)})(1-\x^{(0,1)})}
,
$$
which is the generating function of the vector partition function
\begin{equation}
\label{e:partition}
c(\vec s) =
\# \left\{
       \vec \lambda \in \N^4 \ \Big| \
       \left(
       \begin{matrix}
       1 & 2 & 1 & 0 \\
       1 & 1 & 0 & 1
       \end{matrix}
       \right)
       \vec \lambda = \vec s
   \right\}
.
\end{equation}
This is the same as the example from
\cite[Section~4]{Beck2004fractions}.  For a given $\vec s\in\Z^n$,
the solution set $P_{\vec s}=\left\{\vec\lambda\in\Q^4 \ \big|\
\vec\lambda\ge\vec 0 \text{ and }\left(
       \begin{matrix}
       1 & 2 & 1 & 0 \\
       1 & 1 & 0 & 1
       \end{matrix}
       \right)
       \vec \lambda = \vec s\right\}$ is a two dimensional polytope in
$\Q^4$, so it is helpful to convert it to a full-dimensional
polytope in $\Q^2$ (without changing the number of integer points).
To do this, extend the transformation matrix from~(\ref{e:partition}) to
$$
M =
       \left(
       \begin{matrix}
       1 & 2 & 1 & 0 \\
       1 & 1 & 0 & 1 \\
       0 & 0 & 1 & 0 \\
       0 & 0 & 0 & 1
       \end{matrix}
       \right),
$$
which is unimodular (that is, it has determinant $\pm 1$ and so, as
a linear transformation, it bijectively maps $\Z^4$ to $\Z^4$), and
perform the change of coordinates $\vec \lambda\mapsto
\vec \lambda'=M\vec\lambda$. Then
\begin{align*}
c(\vec s) &= \# \left\{\vec \lambda'\in\Z^4 \ \big| \
M^{-1}\vec \lambda'\ge 0 \text{ and }\left(
       \begin{matrix}
       1 & 2 & 1 & 0 \\
       1 & 1 & 0 & 1
       \end{matrix}
       \right)M^{-1}\vec \lambda'=\vec s\right\}\\
&=\# \left\{\vec \lambda'\in\Z^4 \ \big| \ M^{-1}\vec \lambda'\ge
0
\text{ and }\lambda'_1=s_1,\lambda'_2=s_2\right\}\\
&=\# \Bigg\{
       (\lambda'_3,\lambda'_4) \in \Z^2 \ \ \bigg|\ \
\left(
       \begin{matrix}
-1&2&1&-2\\
1&-1&-1&1\\
0&0&1&0\\
0&0&0&1
       \end{matrix}
       \right)
\left(
\begin{matrix}
s_1\\
s_2\\
\lambda'_3\\
\lambda'_4
\end{matrix}
\right)
       \ge \vec 0
   \Bigg\}
.
\end{align*}
This is the enumeration problem that was our running example in the last section.
\end{example}
\boxend

We will also need the following lemma.

\begin{lemma}
\label{Chambers} Let
$\Phi(m,d)=\binom{m}{0}+\binom{m}{1}+\cdots+\binom{m}{d}$. Then $m$
hyperplanes in $\Q^d$ decompose the space into at most $\Phi(m,d)$
polyhedral chambers.  Furthermore, if we fix $d$, then there is a
polynomial time algorithm which, given $m$ hyperplanes in $\Q^d$,
computes the defining inequalities for each of these chambers.
\end{lemma}

\begin{proof}
This lemma is well known, especially the first part, see, for
example, Section 6.1 of \cite{Matousek02}.  We prove both parts by
induction on $m$. Certainly the statement is true for $m$=0. Suppose
we have a collection of $m$ hyperplanes
$\mathcal{H}_1,\cdots,\mathcal{H}_m$, and assume that these
decompose $\Q^d$ into at most $\Phi(m,d)$ polyhedral chambers whose
defining inequalities may be determined in polynomial time.  Let us
then add a new hyperplane $\mathcal{H}_{m+1}$, which will split some
of the old chambers in two.  The chambers that it splits correspond
exactly to the chambers that the $m$ hyperplanes
$\mathcal{H}_i\cap\mathcal{H}_{m+1}\subset \mathcal{H}_{m+1}$, for
$1\le i\le m$, decompose the $(d-1)$-dimensional space
$\mathcal{H}_{m+1}$ into. Inductively, there are at most
$\Phi(m,d-1)$ of these chambers in $\mathcal{H}_{m+1}$, and their
descriptions may be computed in polynomial time. Therefore, the
hyperplanes $\mathcal{H}_1,\cdots,\mathcal{H}_{m+1}$ decompose
$\Q^m$ into at most $\Phi(m,d)+\Phi(m,d-1)=\Phi(m+1,d)$ chambers,
and we may compute their descriptions in polynomial time.
\end{proof}

A generating function in the form (\ref{GFForm}) is simply the sum
of terms like those in the statement of Corollary \ref{CorVPF}, so the
first half of Theorem \ref{ThmEquivalence} follows from the
following lemma.

\begin{lemma}
Fix $d$.  There is a polynomial time algorithm which, given
piecewise step-polynomials $c_i:\Z^d\rightarrow\Q$, computes $c(\vec
s)=\sum_i c_i(\vec s)$ as a piecewise step-polynomial.
\end{lemma}

\begin{proof}
Suppose $c_i(\vec s)$ are given as piecewise step-polynomials, and
let $c(\vec s)=\sum_i c_i(\vec s)$.  We would like to compute
$c(\vec s)$ as a piecewise step-polynomial.  For each $i$, let
$\{\langle \vec a_{ij}, \vec x\rangle \le b_{ij}\}_j$ be the
collection of linear inequalities that define the chambers of the
piecewise step-polynomial representation of $c_i(\vec s)$. By
Lemma~\ref{Chambers}, we can compute in polynomial time the chambers
in $\Q^n$ determined by the collection of all inequalities
$\{\langle \vec a_{ij}, \vec x\rangle \le b_{ij}\}_{i,j}$.  These
will be the chambers in the piecewise step-polynomial representation
of $c(\vec s)$.  Within a particular chamber, each $c_i(\vec s)$ is
defined by
\[c_i(\vec s)=\sum_{j=1}^{n_i}\alpha_{ij}\prod_{k=1}^{d_{ij}}
\floor{\langle \vec a_{ijk},\vec s\rangle+b_{ijk}},\] where
$\alpha_{ij}\in\Q$, $\vec a_{ijk}\in\Q^d$, and $b_{ijk}\in\Q$, and
so $c(\vec s)=\sum_i c_i(\vec s)$ is simply a sum of such
functions.
\end{proof}

%\label{s:generating}

The first half of Theorem \ref{ThmEquivalence} is now proved. The
second half of Theorem~\ref{ThmEquivalence} depends on
Proposition~\ref{p:generating}, which we prove now.

\begin{proof}[Proof of Proposition~\ref{p:generating}]
\label{p:generating:proof}
Given a parametric polytope $P \subset \Q^{n}\times\Q^{d}$, apply
Theorem 4.4 of \cite{BP99} (see Proposition~\ref{p:core}) directly
on $P$ (that is, not considering $P$ as a parametric polytope but as
a polyhedron in its own right) to obtain the rational generating
function
\[g(P\cap\Z^{n+d};\x,\y)=\sum_{(\vec s,\vec t)\in P\cap\Z^{n+d}}\x^{\vec s}\y^{\vec t}\]
in polynomial time.  Then the generating function $f(\x)$ can be obtained
by substituting $\y=\vec 1$, i.e.,
\[
f(\x) = \sum_{\vec s\in\Z^n}c(\vec s)\x^{\vec s}=
g\big(P\cap\Z^{n+d};\x,\vec 1\big).\] We may perform this
substitution in polynomial time using
Lemma~\ref{l:specialization}. The result is in the form
(\ref{GFForm}).
\end{proof}

Given a piecewise step-polynomial
$c(\vec s)$, we would like to compute the rational generating
function $f(\x)=\sum_{\vec s\in\Z^n}c(\vec s)\x^{\vec s}$.
It suffices to prove it
for functions of the form
\[c(\vec s)=\left\{%
\begin{array}{ll}
    \prod_{j=1}^{d} \floor{\langle \vec a_{j},\vec s\rangle+b_{j}},
    & \hbox{for }\vec s\in Q \\
    0, & \hbox{for }\vec s\notin Q
,
\end{array}%
\right.\] where $Q$ is a rational polyhedron, $\vec a_{j}\in\Q^n$, and
$b_j\in\Q$, because all piecewise step-polynomials may be written as
linear combinations of functions of this form.

Let $P\subset\Q^n\times\Q^d$ be the polyhedron
\[P=\big\{(\vec s, \vec t)\in\Q^n\times\Q^d \ \mid \
\vec s \in Q \text{ and }1\le t_j\le\langle \vec a_{j},\vec
s\rangle+b_{j}, \text{ for }1\le j\le d\big\}.\] Then
\[c(\vec s)=\#\left\{\vec t\in\Z^d \mid (\vec s,\vec t)\in P\right\},\]
and we may compute $f(\x)=\sum_{\vec s}c(\vec s)\x^{\vec s}$ as a
rational generating function using Proposition~\ref{p:generating}.
The second half of Theorem \ref{ThmEquivalence} follows.

Finally, we prove Corollary~\ref{CorParametricProjection}.

\begin{proof}[Proof of Corollary \ref{CorParametricProjection}]
Let
\[S=\left\{(\vec s,\vec t)\in\Z^{n}\times\Z^{d} \mid
\exists \vec u\in\Z^{m}:\
(\vec s,\vec t,\vec u)\in P\right\}.\]
Then we may compute, in polynomial time, the
generating function
\[f(S;\x,\y)=\sum_{(\vec s,\vec t)\in S}\x^{\vec s}\y^{\vec t},\]
using Theorem 1.7 of \cite{BW03}. Next we compute $f\big(S;\x,\vec
1\big)$ using Lemma~\ref{l:specialization}, and the $c(\vec s)$ that
we desire to compute is the piecewise step-polynomial representation
of this generating function. Applying Theorem \ref{ThmEquivalence},
the proof follows (since $P$ is bounded, $\sum_{\vec s}c(\vec s)\x^{\vec s}$ converges everywhere to $f\big(S;\x,\vec
1\big)$, and so any $\vec l$ not orthogonal to any of $\vec b_{ij}$ can be used
in the application of this theorem).
\end{proof}

\section*{Acknowledgements}
We are grateful to Alexander Barvinok for several helpful discussions.  We would also like to thank the anonymous referees for their help in
improving the exposition of this manuscript.

\providecommand{\bysame}{\leavevmode\hbox
to3em{\hrulefill}\thinspace}
\providecommand{\MR}{\relax\ifhmode\unskip\space\fi MR }
\providecommand{\MRhref}[2]{%
  \href{http://www.ams.org/mathscinet-getitem?mr=#1}{#2}
} \providecommand{\href}[2]{#2}

\bibliographystyle{alpha}
\newcommand{\etalchar}[1]{$^{#1}$}

\ \\

\noindent\textsc{Department of Computer Science, K.U. Leuven, Belgium}\\
\emph{Email: }\texttt{Sven.Verdoolaege@cs.kuleuven.be}%
\footnote{Currently at Leiden Institute of Advanced Computer Science,
Universiteit Leiden, The Netherlands, \texttt{sverdool@liacs.nl}}\\

\ \\

\noindent\textsc{Department of Mathematics, Oberlin College, Oberlin, Ohio}\\
\emph{Email: }\texttt{Kevin.Woods@oberlin.edu}
\end{document}